\documentclass[11pt]{article}
\usepackage{amssymb, amsmath}
\usepackage{epsfig}
\usepackage{subfigure}
\usepackage{color}
\usepackage{float}

\newtheorem{theorem}{Theorem}[section]
\newtheorem{definition}[theorem]{Definition}

\newtheorem{example}[theorem]{Example}

\newenvironment{defn*}{\begin{definition}}{\end{definition}}

\begin{document}

%DO NOT EDIT ABOVE THIS LINE EXCEPT FOR THE NAMES AND THE DATE%

\title{\bf On the invariant manifolds of the fixed point of a second order nonlinear difference equation}

\author{Mehmet TURAN}
\date{}
\maketitle

\begin{center}
{\it Atilim University, Department of
Mathematics,  Incek  06836, Ankara, Turkey}\\
{\it e-mail: mehmet.turan@atilim.edu.tr}\\
{\it Tel: +90 312 586 8585,  Fax: +90 312 586 8091}
\end{center}

\begin{abstract} 

This paper addresses the asymptotic approximations of the stable and
unstable manifolds for the saddle fixed point and the 2-periodic
solutions of the difference equation $x_{n+1} = \alpha + \beta x_{n-1}+x_{n-1}/x_{n},$
where $\alpha>0,$ $0\leqslant \beta <1$ and the initial conditions $x_{-1}$
and $x_0$ are positive numbers. These manifolds determine
completely global dynamics of this equation. The theoretical results are supported by some numerical examples.


\end{abstract}

\section{Introduction}

Many real world processes are studied by means of difference equations. Because
of their wide range of applications in mechanics, economics, electronics, chemistry,
ecology, biology, etc., the theory of discrete dynamical systems has been under intensive
development and many researchers have been paying their attention to the study
of these systems \cite{agarwal, elaydi, HK, kelley, kocic, ladas, KM, laks2, hassan, Wiggins}.

The equation
\begin{align}
x_{n+1}=\alpha + \frac{x_{n-1}}{x_n}, \qquad n=0,1,2,\ldots \label{eqn}
\end{align}
was investigated by many researchers. 
The equation \eqref{eqn}, for $\alpha \in [0,\infty)$ 
and the initial conditions $x_{-1}$ and $x_0$ 
being arbitrary positive real numbers, has been considered in \cite{amleh}. 
There, the authors analyzed the global stability, the
boundedness character, and the periodic nature of the positive solutions of \eqref{eqn}.
The global stability, the permanence, 
and the oscillation character of the recursive equation \eqref{eqn} 
for nonnegative values of the parameter $\alpha$ 
with negative initial conditions $x_{-1}$ and $x_0$ was investigated in \cite{hamza}. 
The same equation for $\alpha< 0$ was taken into account in \cite{stevic2}.
The global bifurcation result for \eqref{eqn} was obtained in \cite{burgic} and 
the asymptotic approximations of the stable and unstable manifolds of the fixed point of 
\eqref{eqn} were discussed in \cite{kul}.

In this work, we consider the difference equation
\begin{align}
x_{n+1}=\alpha+\beta x_{n-1} + \frac{x_{n-1}}{x_n}, \qquad n=0,1,2,\ldots \label{main}
\end{align}
where $\alpha \geqslant 0,$ $0 \leqslant \beta < 1,$ 
and the initial conditions $x_{-1}$ and $x_0$ are positive real numbers. 
Clearly, when $\beta=0,$ the equation \eqref{main} reduces to \eqref{eqn}. For this reason,
the results obtained in the current paper covers those given in \cite{kul}.

Equation \eqref{main}
has the unique fixed point 
\begin{align} \label{fp}
\bar{x}=\frac{1+\alpha}{1-\beta}.
\end{align}

Letting $y_n=x_{n-1}$ and $z_n=x_n,$ \eqref{main}
can be written as
\begin{align} \label{mainsys}
\begin{array}{l}
y_{n+1} = z_n \\
z_{n+1} = \alpha + \beta y_n + \displaystyle\frac{y_n}{z_n}
\end{array}
\end{align}
together with the initial conditions $y_0=x_{-1},$ $z_0=x_0.$ Introducing the mapping
\begin{align} \label{t}
T\begin{pmatrix}
y \\ z
\end{pmatrix} =
\begin{pmatrix}
z \\ \alpha + \beta y + \displaystyle\frac{y}{z}
\end{pmatrix},
\end{align}
\eqref{mainsys} is written as 
\begin{align*}
\begin{pmatrix}
y_{n+1} \\ z_{n+1}
\end{pmatrix}
=
T\begin{pmatrix}
y_n \\ z_n
\end{pmatrix}.
\end{align*}
$T$ has a unique fixed point $(\bar{x}, \bar{x})$ where $\bar{x}$ is given by \eqref{fp}.

The following result for \eqref{main} was given in \cite{at}:
\begin{theorem} Let $0\leqslant \beta < 1.$ For the equation \eqref{main}, one has:
\begin{itemize} 
\item If $0\leqslant \alpha < 1,$ the equilibrium point $\bar{x}$ is unstable;
\item If $\alpha = 1,$  then there exists periodic solutions with period 2. Moreover, any non periodic solution of \eqref{main} converges either to the fixed point or to a two-periodic solution;
\item If $\alpha >1,$ then the equilibrium point $\bar{x}$ is globally asymptotically stable.
\end{itemize}
\end{theorem}
To complete the global dynamics of \eqref{main}, the present paper addresses the equations of stable and unstable manifolds of
the equilibrium solution and the stable manifold of period-two solutions of \eqref{main}.
The following definition of the stable and unstable manifolds and the next theorem about their existence can be found in \cite[Definition 15.18, Theorem 15.19, pp.457]{HK} and also in \cite{mars}. We present these only with a minor change in notations for the convenience of the present paper.

\begin{defn*} Let ${\mathcal N}$ be a neighborhood 
of a fixed point $\bar{x}$ of a diffeomorphism $T$ defined in ${\mathcal N}.$ 
Then, the local stable manifold $W^s(\bar{x}, {\mathcal N}),$ and
the local unstable manifold $W^u(\bar{x}, {\mathcal N})$ of $\bar{x}$ are defined, 
respectively, to be the following subsets of ${\mathcal N}:$ 
\begin{align*}
W^s(\bar{x}, {\mathcal N})=\{{\bf x}\in {\mathcal N}: T^n({\bf x})\in {\mathcal N}, \ \text{for all } \ n\geqslant 0,
\ \text{and} \ T^n({\bf x})\to \bar{x}, \ \text{as} \ n\to\infty \} \\
W^u(\bar{x}, {\mathcal N})=\{{\bf x}\in {\mathcal N}: T^{-n}({\bf x})\in {\mathcal N}, \ \text{for all } \ n\geqslant 0,
\ \text{and} \ T^{-n}({\bf x})\to \bar{x}, \ \text{as} \ n\to\infty \}
\end{align*}
\end{defn*}

\begin{theorem}[Stable and Unstable Manifolds]\label{thminv}
Let $T$ be be a diffeomorphism with a hyperbolic saddle point $\bar{x},$ that is, the linearized
map $DT(\bar{x})$ at the fixed point has nonzero eigenvalues $|\lambda_1|<1$ and $|\lambda_2|>1.$
Then $W^s(\bar{x}, {\mathcal N})$ is a curve tangent at $\bar{x}$ to, and a graph over, the eigenspace
corresponding to $\lambda_1,$ while $W^u(\bar{x}, {\mathcal N})$ is a curve tangent at $\bar{x}$ to, and a graph
over, the eigenspace corresponding to $\lambda_2$. These curves are as smooth as
the map $T.$
\end{theorem}
For the following theorem see \cite[Theorem 6, pp 34]{carr}.
\begin{theorem}[Center Manifold]\label{thmc}
Let $T:{\mathbb R}^{n+m}\to{\mathbb R}^{n+m}$ have the following form:
$$T(x,y)=(Ax+f(x,y), \ Bx+g(x,y))$$
where $x\in{\mathbb R}^{n},$ $y\in{\mathbb R}^{m},$ $A$ and $B$ are square matrices such that each eigenvalue of $A$ has modulus 1 and each eigenvalue of $B$ has modulus less than 1, $f$ and $g$ are $C^2$ and $f, g$ and their first order derivatives are zero at the origin. Then, there exists a center manifold $h:{\mathbb R}^n\to{\mathbb R}^m$ for $T.$
More precisely, for some $\varepsilon >0$ there exists a $C^2$ function $h:{\mathbb R}^n\to{\mathbb R}^m$ with 
$h(0)=h'(0)=0$ such that $|x|<\varepsilon$ and $(x_1, y_1)=T(x,h(x))$ implies $y_1=h(x_1).$
\end{theorem}

The paper is organized as follows: In the next chapter, the normal form of the map $T$ and the equations of unstable and stable manifolds of the equilibrium solution are given. Chapter 3 deals with the normal form and invariant manifolds of the map $T^2.$ Finally, Chapter 4 is devoted to some numerical examples to illustrate the theoretical results.

\section{Normal form and invariant manifolds of the map $T$}

\subsection{Normal Form}

To obtain the normal form of the map $T,$ first, we transform its fixed point to the origin. For this, let $u_n=y_n-\bar{x}$ and $v_n=z_n-\bar{x}.$ Then, \eqref{mainsys} becomes
\begin{align} \label{mainsys2}
\begin{array}{l}
u_{n+1}=v_n\\
v_{n+1}=\beta u_n+\displaystyle\frac{u_n-v_n}{v_n+\bar{x}}
\end{array}.
\end{align}
For the mapping
\begin{align*}
F\begin{pmatrix}
u \\ v
\end{pmatrix}
=
\begin{pmatrix}
v\\ \beta u+ \frac{u-v}{v+\bar{x}}
\end{pmatrix},
\end{align*}
\eqref{mainsys2} is written as 
\begin{align*}
\begin{pmatrix}
u_{n+1} \\ v_{n+1}
\end{pmatrix}
=
F\begin{pmatrix}
u_n \\ v_n
\end{pmatrix}.
\end{align*}
The Jacobian matrix of $F$ at its unique fixed point $(0,0)$ is
\begin{align*}
J=
\begin{pmatrix}
0 & 1\\ \beta+\frac{1}{\bar{x}} & -\frac{1}{\bar{x}}
\end{pmatrix}
\end{align*}
which has the eigenvalues 
\begin{align}
\lambda_{1} = \frac{-1-\theta}{2\bar{x}} 
\quad \text{and} \quad 
\lambda_{2} = \frac{-1+\theta}{2\bar{x}}\label{evalues}
\end{align}
with the corresponding eigenvectors
\begin{align}
{\bf v}_1 = \left(\frac{-2\bar{x}}{1+\theta} \:, \: \: 1\right)^T  
\quad \text{and} \quad 
{\bf v}_2 = \left(\frac{-2\bar{x}}{1-\theta} \:, \: \: 1\right)^T  
\label{evectors}
\end{align}
respectively, where $\theta=\sqrt{1+4\bar{x}+4\beta\bar{x}^2}.$
Thus,
\begin{align*}
F
\begin{pmatrix}
u \\ v
\end{pmatrix}
= J\cdot
\begin{pmatrix}
u \\ v
\end{pmatrix}
+ H
\begin{pmatrix}
u \\ v
\end{pmatrix}
\end{align*}
where
\begin{align*}
H
\begin{pmatrix}
u \\ v
\end{pmatrix}
=
\begin{pmatrix}
0 \\ \frac{v(v-u)}{\bar{x}(v+\bar{x})}
\end{pmatrix}.
\end{align*}
Thus, \eqref{mainsys} is equivalent to
\begin{align} \label{maineq}
\begin{pmatrix}
u_{n+1}\\ v_{n+1} 
\end{pmatrix}
=
\begin{pmatrix}
0 & 1 \\ \beta+\frac{1}{\bar{x}} & -\frac{1}{\bar{x}}
\end{pmatrix}
\begin{pmatrix}
u_n\\ v_n 
\end{pmatrix}
+H
\begin{pmatrix}
u_n\\ v_n 
\end{pmatrix}.
\end{align}
Set 
$P=({\bf v}_1 \: {\bf v}_2),$ where ${\bf v}_1$ and ${\bf v}_2$ are given by \eqref{evectors}, and let 
\begin{align} \label{uvp}
\begin{pmatrix}
u_n \\ v_n
\end{pmatrix} = P\cdot \begin{pmatrix}
\xi_n \\ \eta_n
\end{pmatrix}.
\end{align}
Then, \eqref{maineq} becomes
\begin{align} 
\begin{pmatrix}
\xi_{n+1} \\ \eta_{n+1}
\end{pmatrix} 
=
\begin{pmatrix}
\lambda_1 & 0 \\ 0 & \lambda_2
\end{pmatrix}
\begin{pmatrix}
\xi_n \\ \eta_n
\end{pmatrix}
+
\begin{pmatrix}
f(\xi_n, \eta_n) \\ g(\xi_n, \eta_n)
\end{pmatrix}
\label{normal}
\end{align}
where
\begin{align} \label{fg}
\begin{array}{l}
f(\xi, \eta)=
\displaystyle\frac{(1+2\beta\bar{x})(\xi+\eta)^2+\theta(\xi^2-\eta^2)}{\theta(\theta-1)(\xi+\eta+\bar{x})},\\
g(\xi, \eta)=
\displaystyle\frac{(1+2\beta\bar{x})(\xi+\eta)^2+\theta(\xi^2-\eta^2)}{\theta(\theta+1)(\xi+\eta+\bar{x})}.
\end{array}
\end{align}
System \eqref{normal} is called the normal form of \eqref{mainsys}.

\subsection{Unstable manifold of the equilibrium solution}

Let $0<\alpha<1.$ Then, as it is stated in \cite{at}, the fixed point $\bar{x}$ of \eqref{main} is unstable. In fact, it can be shown that $|\lambda_1|>1$ and $|\lambda_2|<1.$ Then, by Theorem \ref{thminv}, there is an unstable manifold $W^u$ which is the graph of an analytic map $\varphi:E_1 \to E_2$
such that $\varphi(0)=\varphi'(0)=0.$ Let $$\varphi(\xi)=a_2\xi^2+a_3\xi^3+O(\xi^4), \quad a_2, a_3 \in{\mathbb R}.$$ Now, we shall compute $a_2$ and $a_3.$ On the manifold $W^u,$  we have $\eta_n=\varphi(\xi_n)$ for $n\in {\mathbb N}_0.$ Thus, the function $\varphi$ must satisfy 
\begin{align}
\varphi(\lambda_1\xi+f(\xi,\varphi(\xi)))=\lambda_2\varphi(\xi)+g(\xi, \varphi(\xi)) \label{unsman}
\end{align}
where $f$ and $g$ are given in \eqref{fg}. Rewriting \eqref{unsman} as a polynomial equation in $\xi$ and equation the coefficients of $\xi^2$ and $\xi^3$ to 0, we obtain
\begin{align}
a_2=\frac{1+\theta+2\beta\bar{x}}{\theta(\theta+1)(\lambda_1^2-\lambda_2)\bar{x}} \label{a2}
\end{align}
and
\begin{align}
a_3=\frac{a_2}{(\lambda_1^3-\lambda_2)\bar{x}}\left[\lambda_2-\lambda_1^2-\frac{1+2\beta\bar{x}}{\theta\bar{x}}\left(\frac{1}{\lambda_1}+\frac{\lambda_1}{\lambda_2}\right)-\frac{\lambda_1}{\lambda_2\bar{x}}\right]. \label{a3}
\end{align}

The local unstable manifold is obtained locally as the graph of the map
$\varphi(\xi)=a_2\xi^2+a_3\xi^3.$ Since $\eta_n=a_2\xi_n^2+a_3\xi_n^3,$ using \eqref{uvp} and $u_n=x_{n-1}-\bar{x},$ $v_n=x_n-\bar{x},$ we can approximate locally the local unstable manifold $W_{loc}^u$ of \eqref{main} as the graph of $\tilde{\varphi}(x)$ such that $U(x, \tilde{\varphi}(x))=0$ where
\begin{align}
U(x,y):=\gamma_1(x-\bar{x})-\gamma_2(y-\bar{x})
+a_2\left[\gamma_1(x-\bar{x})+\gamma_3(y-\bar{x})\right]^2
+a_3\left[\gamma_1(x-\bar{x})+\gamma_3(y-\bar{x})\right]^3 \label{manU}
\end{align}
in which 
\begin{align} \label{gam}
\gamma_1=\frac{1+\beta\bar{x}}{-\theta},\quad 
\gamma_2=\frac{\theta-1}{2\theta} \quad \text{and} \quad 
\gamma_3=\frac{\theta+1}{2\theta}.
\end{align}

It is easy to see that the function $\tilde{\varphi}(x)$ satisfies 
\begin{align*}
\tilde{\varphi}(\bar{x})=\bar{x} \quad \textnormal {and} \quad 
\tilde{\varphi}'(\bar{x})=-\frac{1+\theta}{2\bar{x}}.
\end{align*}
Thus, we have proved the following theorem:
\begin{theorem} \label{thmusman}
The local unstable manifold of \eqref{main} corresponding to the saddle point $\bar{x}$ has the asymptotic equation $U(x,\tilde{\varphi}(x))=0$ where $U(x,y)$ is given by \eqref{manU}.
\end{theorem}

\subsection{Stable manifold of the equilibrium solution}

Since $|\lambda_1|>1$ and $|\lambda_2|<1,$ by Theorem \ref{thminv}, there is a stable manifold $W^s$ which is the graph of an analytic map $\psi:E_1 \to E_2$
such that $\psi(0)=\psi'(0)=0.$ Let $$\psi(\eta)=b_2\eta^2+b_3\eta^3+O(\eta^4), \quad b_2, b_3 \in{\mathbb R}.$$ Now, we shall compute the coefficients $b_2$ and $b_3.$ On the manifold $W^s,$  we have $\xi_n=\psi(\eta_n)$ for $n\in {\mathbb N}_0.$ Thus, the function $\psi$ must satisfy 
\begin{align}
\psi(\lambda_2\eta+g(\psi(\eta),\eta))=\lambda_1\psi(\eta)+f(\psi(\eta),\eta) \label{sman}
\end{align}
where $f$ and $g$ are given in \eqref{fg}. Rewriting \eqref{sman} as a polynomial equation in $\eta$ and equating the coefficients of $\eta^2$ and $\eta^3$ to 0, we obtain
\begin{align}
b_2=\frac{1-\theta+2\beta\bar{x}}{\theta(\theta-1)(\lambda_2^2-\lambda_1)\bar{x}} \label{b2}
\end{align}
and
\begin{align}
b_3=\frac{b_2}{(\lambda_b^3-\lambda_1)\bar{x}}\left[\lambda_1-\lambda_2^2+\frac{1+2\beta\bar{x}}{\theta\bar{x}}\left(\frac{1}{\lambda_2}+\frac{\lambda_2}{\lambda_1}\right)-\frac{\lambda_2}{\lambda_1\bar{x}}\right]. \label{b3}
\end{align}

The local stable manifold is obtained locally as the graph of the map
$\psi(\eta)=b_2\eta^2+b_3\eta^3.$ Since $\xi_n=b_2\eta_n^2+b_3\eta_n^3,$ using \eqref{uvp} and $u_n=x_{n-1}-\bar{x},$ $v_n=x_n-\bar{x},$ we can approximate locally the local stable manifold $W_{loc}^s$ of \eqref{main} as the graph of $\tilde{\psi}(y)$ such that $S(\tilde{\psi}(y), y)=0$ where
\begin{align}
S(x,y):=\gamma_1(x-\bar{x})+\gamma_3(y-\bar{x})
-b_2\left[\gamma_1(x-\bar{x})-\gamma_2(y-\bar{x})\right]^2
+b_3\left[\gamma_1(x-\bar{x})-\gamma_2(y-\bar{x})\right]^3. \label{manS}
\end{align}
It is easy to see that the function $\tilde{\psi}(x)$ satisfies 
\begin{align*}
\tilde{\psi}(\bar{x})=\bar{x} \quad \textnormal {and} \quad 
\tilde{\psi}'(\bar{x})=\frac{2\bar{x}}{\theta-1}.
\end{align*}
Thus, we have proved the following theorem:
\begin{theorem} \label{thmsman}
The local stable manifold of \eqref{main} corresponding to the saddle point $\bar{x}$ has the asymptotic equation $S(\tilde{\psi}(y),y)=0$ where $S(x,y)$ is given by \eqref{manS}.
\end{theorem}

\section{Normal form and invariant manifold of the map $T^2$}

\subsection{Normal Form}

For the map $T$ given by \eqref{t}, one has
\begin{align}\label{t2}
T^2\begin{pmatrix}
y \\ z
\end{pmatrix}
=
\begin{pmatrix}
1+\beta y + y/z \\
1+\beta z + \frac{z}{1+\beta y + y/z}
\end{pmatrix}.
\end{align}
That is,
\begin{align}
\begin{pmatrix}
y_{n+2} \\ z_{n+2}
\end{pmatrix}
=
T^2\begin{pmatrix}
y_n \\ z_n
\end{pmatrix}. \label{yzt2}
\end{align}

Firstly, we note that the fixed point $(\bar{x},\bar{x})$ of $T$ is also a fixed point of $T^2.$ That is why, in this section, the fixed point $(\bar{x},\bar{x})$ of $T^2$ is ignored, and the main focus will be on the other fixed points. As it was shown in \cite{at}, when $\alpha=1,$ equation \eqref{main} has infinitely many 2-periodic solutions each of which corresponds to a fixed point of $T^2.$ Indeed, if $\Phi>1/(1-\beta),$ then all the fixed points of $T^2$ are given by $(\Phi, \Psi)$ where $\Psi=\Phi/[(1-\beta)\Phi-1]).$ It is worth mentioning that $\Psi>1/(1-\beta)$ and for the initial conditions $x_{-1}=\Phi$ and $x_0=\Psi,$ the solution of \eqref{main} is $\{\Phi, \Psi, \Phi, \Psi, \ldots\}.$ Swapping the initial values produces the periodic solution $\{\Psi, \Phi, \Psi, \Phi, \ldots\}.$

  As it was done in the previous section, the fixed point $(\Phi, \Psi)$ will be transformed to the origin. For this, let $u=y-\Phi$ and $v=z-\Psi.$ Then, we get the map
\begin{align} \label{F0}
F_0\begin{pmatrix}
u \\ v
\end{pmatrix}
:=T^2\begin{pmatrix}
u+\Phi \\ v+\Psi
\end{pmatrix}
-\begin{pmatrix}
\Phi \\ \Psi
\end{pmatrix}
=\begin{pmatrix}
\beta u + \frac{u + \Phi}{v+\Psi} - \frac{\Phi}{\Psi}\\
\beta v + \frac{(v + \Psi)^2}{v+\Psi+(u+\Phi)(1+\beta v+\beta\Psi)} - \frac{\Psi}{\Phi}
\end{pmatrix},
\end{align}
for which \eqref{yzt2} can be written as
\begin{align}
\begin{pmatrix}
u_{n+2} \\ v_{n+2}
\end{pmatrix}
=
F_0\begin{pmatrix}
u_n \\ v_n
\end{pmatrix}. \label{unvn}
\end{align}

It is clear to see that $(0,0)$ is a fixed point of $F_0$ and the Jacobian of $F_0$ at this fixed point is
\begin{align}\label{J0}
J_0=\begin{pmatrix}
\beta+\frac{1}{\Psi} & -\frac{\Phi}{\Psi^2} \\
-\frac{\beta \Psi+1}{\Phi^2} & \beta+\frac{1}{\Phi}+\frac{1}{\Psi\Phi}
\end{pmatrix}.
\end{align}
Using the relation $\frac{1}{\Phi}+\frac{1}{\Psi}=1-\beta,$ one can rewrite \eqref{J0} as
\begin{align*}
J_0=\begin{pmatrix}
1-\frac{1}{\Phi} & -\frac{\Phi}{\Psi^2} \\
\frac{\Psi(1-\Phi)}{\Phi^3} & 1-\frac{1}{\Psi}+\frac{1}{\Psi\Phi}
\end{pmatrix}.
\end{align*}
from which it can be derived that the eigenvalues of $J_0$ are 
\begin{align*}
\lambda_{01}=\left(1-\frac{1}{\Phi}\right)\left(1-\frac{1}{\Psi}\right) \quad \text{and} \quad
\lambda_{02}=1.
\end{align*}
Since $0<\lambda_{01}<1$ and $\lambda_{02}=1,$ the fixed point $(\Phi, \Psi)$ is stable. The eigenvectors corresponding to the eigenvalues $\lambda_{01}$ and $\lambda_{02}$ are 
\begin{align}
{\bf v}_{01}= \left(\frac{\Phi^2}{(\Phi-1)\Psi} \:, \: \: 1\right)^T \quad \text{and} \quad
{\bf v}_{02}= \left(-\frac{\Phi^2}{\Psi^2} \:, \:\: 1\right)^T, \label{ev0}
\end{align}
respectively. Thus, the map obtained in \eqref{F0} can be written as
\begin{align} \label{F02}
F_0\begin{pmatrix}
u \\ v
\end{pmatrix}
=J_0 \cdot 
\begin{pmatrix}
u \\ v
\end{pmatrix}
+H_0\begin{pmatrix}
u \\ v
\end{pmatrix}
\end{align}
where
\begin{align*}
H_0\begin{pmatrix}
u \\ v
\end{pmatrix}
=\begin{pmatrix}
\frac{v(\Phi v-\Psi u)}{\Psi^2(v+\Psi)} \\ 
\frac{(v + \Psi)^2}{v+\Psi+(u+\Phi)(1+\beta v+\beta\Psi)} + \frac{(\beta\Psi+1)u}{\Phi^2} - \frac{(\Psi+1)v}{\Psi\Phi} - \frac{\Psi}{\Phi} 
\end{pmatrix}.
\end{align*}
Therefore, \eqref{yzt2} is equivalent to
\begin{align}
\begin{pmatrix}
u_{n+2} \\ v_{n+2}
\end{pmatrix}
=
\begin{pmatrix}
1-\frac{1}{\Phi} & -\frac{\Phi}{\Psi^2} \\
\frac{\Psi(1-\Phi)}{\Phi^3} & 1-\frac{1}{\Psi}+\frac{1}{\Psi\Phi}
\end{pmatrix}
\begin{pmatrix}
u_n \\ v_n
\end{pmatrix}
+ 
H_0 \begin{pmatrix}
u_n \\ v_n
\end{pmatrix}. \label{unvn2}
\end{align}
Set 
$P_0=({\bf v}_{01} \: {\bf v}_{02}),$ where ${\bf v}_{01}$ and ${\bf v}_{02}$ are given by \eqref{ev0}, and let 
\begin{align} \label{uvp0}
\begin{pmatrix}
u \\ v
\end{pmatrix} = P_0\cdot \begin{pmatrix}
\xi \\ \eta
\end{pmatrix}.
\end{align} 
Then, \eqref{unvn2} leads to
\begin{align} 
\begin{pmatrix}
\xi_{n+2} \\ \eta_{n+2}
\end{pmatrix} 
=
\begin{pmatrix}
\lambda_{01} & 0 \\ 0 & \lambda_{02}
\end{pmatrix}
\begin{pmatrix}
\xi_n \\ \eta_n
\end{pmatrix}
+
\begin{pmatrix}
f_0(\xi_n, \eta_n) \\ g_0(\xi_n, \eta_n)
\end{pmatrix} \label{norm0}
\end{align}
where
\begin{align} \label{f0g0}
\begin{array}{l}
f_0(\xi, \eta)=
\displaystyle \frac{\Phi-1}{\Psi+\Phi-1}\left( \zeta-\frac{\xi}{\Phi\Psi}-\frac{\eta(\Phi+\Psi)}{\Phi(\xi+\eta+\Psi)}+\frac{(\xi+\eta)\xi}{\Phi(1-\Phi)(\xi+\eta+\Psi)}\right), \\ 
g_0(\xi, \eta)=
\displaystyle \frac{\Psi}{\Psi+\Phi-1}\left(\zeta-\frac{\eta(\Phi+\Psi)}{\Phi\Psi}-\frac{\xi}{\Phi(\xi+\eta+\Psi)}+\frac{(1-\Phi)(\Phi+\Psi)(\xi+\eta)\eta}{\Phi\Psi^2(\xi+\eta+\Psi)}\right).
\end{array}
\end{align}
and 
$$\zeta=\frac{(\xi+\eta+\Psi)^2}{\xi+\eta+\Psi+\left(\frac{\Phi^2\xi}{(\Phi-1)\Psi}-\frac{\Phi^2\eta}{\Psi^2}+\Phi\right)(1+\beta \xi+\beta\eta+\beta\Psi)}-\frac{\Psi}{\Phi}.$$
System \eqref{norm0} is the normal of \eqref{yzt2}.

\subsection{Stable set of 2-periodic solution $\{(\Phi, \Psi), (\Psi, \Phi)\}$}

Since $0<\lambda_{01}<1$ and $\lambda_{02}=1,$ by Theorem \ref{thmc}, there is an invariant curve ${\mathcal C}$ (called center manifold) which is the graph of an analytic map $h$ such that $h(0)=h'(0)=0.$ Let 
\begin{align*}
h(\xi)=c_2\xi^2+c_3\xi^3+O(\xi^4), \quad c_2, c_3 \in{\mathbb R}.
\end{align*}
Now, we shall compute $c_2$ and $c_3.$ The function $h$ must satisfy 
\begin{align}
h(\lambda_{01}\xi+f_0(\xi, h(\xi)))=h(\xi)+g_0(\xi,h(\xi)) \label{cman}
\end{align} 
where $\lambda_{01}$ is given in \eqref{ev0}, $f_0$ and $g_0$ are given in \eqref{f0g0}. Rewriting \eqref{cman} as a polynomial equation in $\xi$ and equation the coefficients of $\xi^2$ and $\xi^3$ to 0, we obtain
\begin{align}\label{c2}
c_2=\frac{\Phi}{(1-\Phi)(\Phi+\Psi-1)(2\Phi\Psi-\Phi-\Psi+1)}
\end{align}
and
\begin{align}\label{c3}
c_3=\frac{\Phi^4\kappa_1+\Phi^3\kappa_2-\Phi^2\kappa_3+\kappa_4\Phi}
{\Psi(\Psi+\Phi-1)(1-\Phi)\kappa_5}\, c_2
\end{align}
where
\begin{align*}
\kappa_1&=3\Psi^2-4\Psi+1\\
\kappa_2&=3\Psi^3-12\Psi^2+9\Psi-1\\
\kappa_3&=5\Psi^3-14\Psi^2+6\Psi+1\\
\kappa_4&=2\Psi^3-4\Psi^2+\Psi+1\\
\kappa_5&=\Phi^2(3\Psi^2-3\Psi+1)-\Phi(3\Psi^2-5\Psi+2)+(\Psi-1)^2.
\end{align*}
Using $\eta_{2n}=c_2\xi_{2n}^2+c_3\xi_{2n}^3,$ the relation \eqref{uvp0} together with $u_{2n}=x_{2n-2}-\Phi$ and  $v_{2n}=x_{2n}-\Psi,$ we can approximate locally the invariant curve ${\mathcal C}$ as the graph of $\tilde{h}(x)$ such that $C(x, \tilde{h}(x))=0$ where
\begin{align*}
C(x,y;\Phi):=\delta_1(x-\Phi)-\delta_2(y-\Psi)
+c_2\left[\delta_1(x-\Phi)+\delta_3(y-\Psi)\right]^2
+c_3\left[\delta_1(x-\Phi)+\delta_3(y-\Psi)\right]^3,
\end{align*}
in which 
\begin{align} \label{del}
\delta_1=\frac{\Psi^2(\Phi-1)}{\Phi^2(\Phi+\Psi-1)},\quad 
\delta_2=\frac{\Psi}{\Phi+\Psi-1} \quad \text{and} \quad 
\delta_3=\frac{\Phi-1}{\Phi+\Psi-1}.
\end{align}

It is easy to see that the function $\tilde{h}(x)$ satisfies 
\begin{align*}
\tilde{h}(\Phi)=\Psi \quad \textnormal {and} \quad 
\tilde{h}'(\Phi)=\frac{\Psi(\Phi-1)}{\Phi^2}.
\end{align*}
Thus, we have proved the following theorem:
\begin{theorem} \label{thmcman}
Let $\Phi>1/(1-\beta)$ and $\Psi=\Phi/[(1-\beta)\Phi-1].$ Then corresponding to the non-hyperbolic period-two solution $\{(\Phi, \Psi), (\Psi, \Phi)\},$ there is an invariant curve which is the union of two curves that are locally given with the asymptotic expansions $C(x,\tilde{h}(x);\Phi)=0$ and $C(x,\tilde{h}(x);\Psi)=0.$
\end{theorem}

\section{Numerical Examples}

In this section, some illustrative examples supporting the theoretical results presented in this article will be constructed. To compare the current results with those given in \cite{kul}, we first take the parameter values as in \cite{kul}. 

\begin{example}
For $\alpha=0.2,$ $\beta=0,$ which is the case $p=0.2$ in \cite[Section 3.4]{kul}, we have
\begin{align*}
U_1(x,y)&= -0.4152273992 x + 0.8491364395 - 0.2923863004 y \\
& \qquad  +0.2419777563(-0.4152273992 x - 0.3508635604 + 0.7076136995 y)^2\\
& \qquad  -0.0974600586(-0.4152273992 x - 0.3508635604 + 0.7076136995 y)^3,\\
S_1(x,y)&=-0.4152273992 x - 0.3508635604 + 0.7076136995 y \\
& \qquad  +0.1961061968(-0.4152273992 x + 0.8491364395 - 0.2923863004 y)^2\\
& \qquad  +0.09806508071(-0.4152273992 x + 0.8491364395 - 0.2923863004 y)^3.
\end{align*}
and 
\begin{align*}
U_2(x,y)&= -0.3492151478 x + 1.214293633 -0.3253924261 y \\ 
& \qquad +0.3059452562(-0.3492151478 x - 0.5857063670 + 0.6746075740y)^2 \\
& \qquad -0.1066716833(-0.3492151478 x - 0.5857063670 + 0.6746075740y)^3, \\
S_2(x,y)&= -0.3492151478 x - 0.5857063670 + 0.6746075740 y \\
& \qquad + 0.1446549340(-0.3492151478 x + 1.214293633 - 0.3253924261y)^2 \\
& \qquad + 0.0525187072(-0.3492151478 x + 1.214293633 - 0.3253924261y)^3.
\end{align*}

Figure \ref{fig1} shows the graphs of the functions $U_1(x,y)=0,$ $S_1(x,y)=0,$ $U_2(x,y)=0$ and $S_2(x,y)=0$ together with a typical trajectory. As it can be seen, the trajectory follows the unstable manifold in both cases.

\begin{figure}[!ht]
\begin{center}
\subfigure[Graphs of $U_1(x,y)=0$ (blue) and $S_1(x,y)=0$ (red) for $\alpha=0.2, \beta=0$]{
\resizebox*{60mm}{!}{\includegraphics{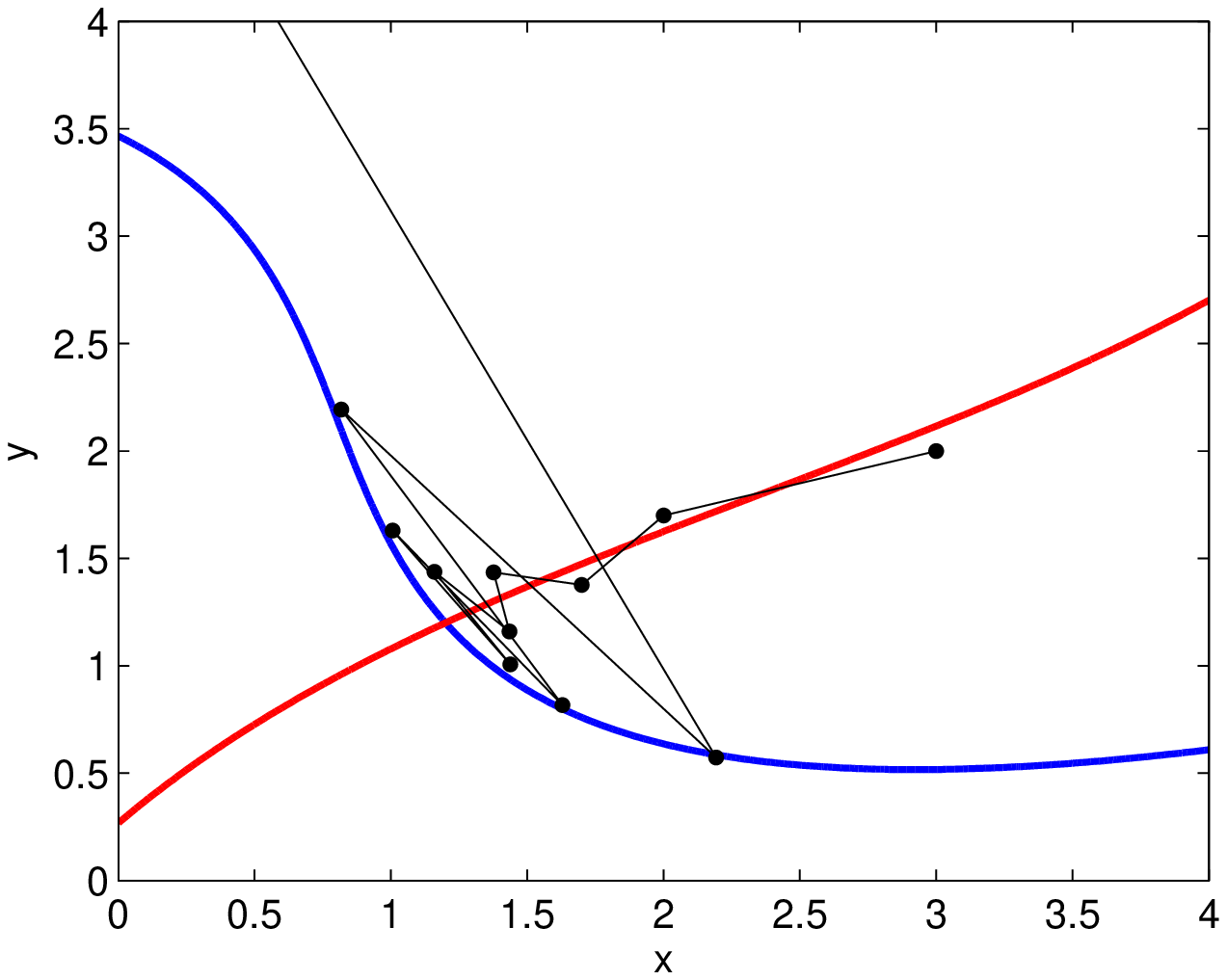}}}\hspace{8mm}%
\subfigure[Graphs of $U_2(x,y)=0$ (blue) and $S_2(x,y)=0$ (red) for $\alpha=0.8, \beta=0$]{
\resizebox*{60mm}{!}{\includegraphics{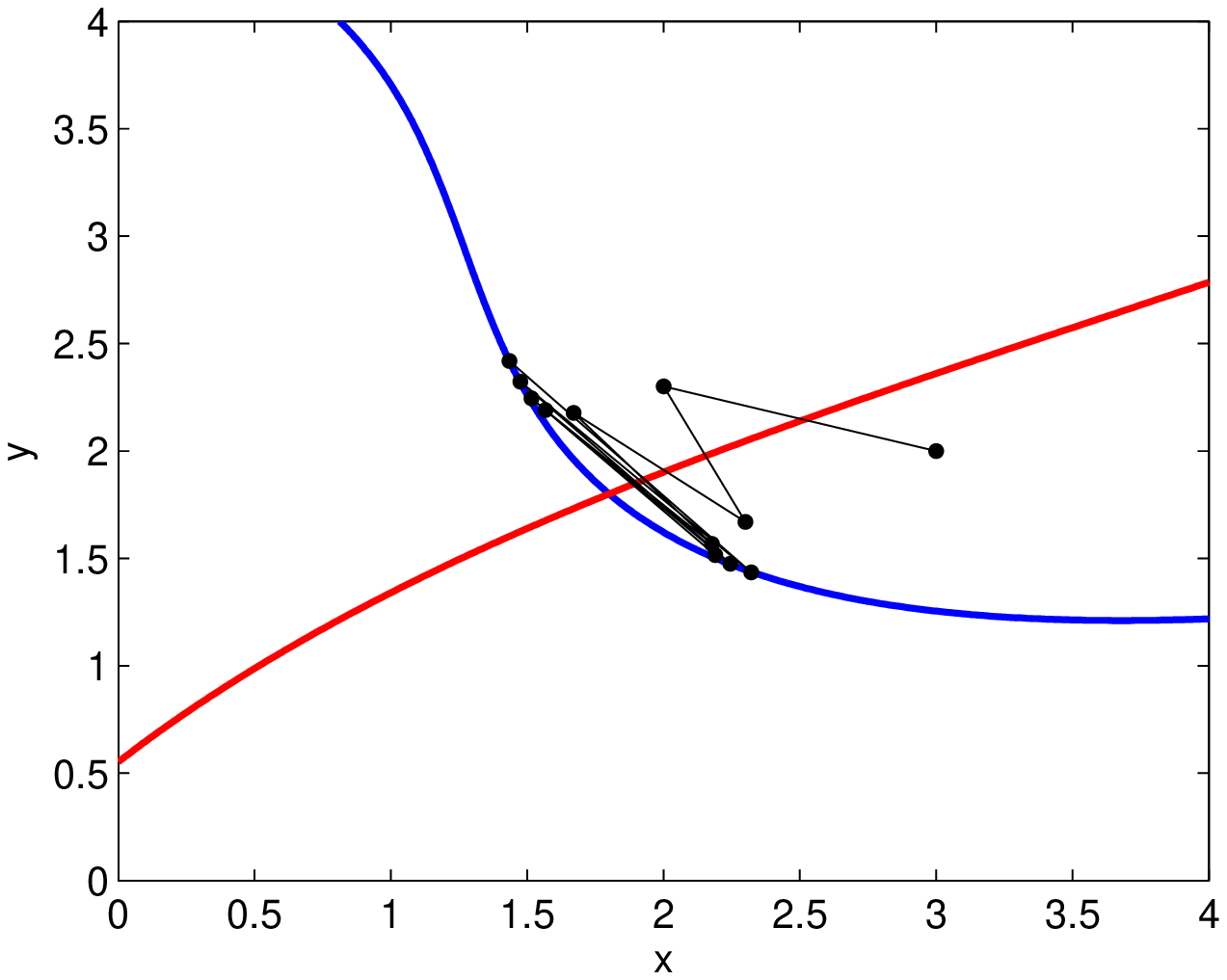}}}
\caption{\label{fig1} Graphs of stable and unstable manifolds together with a typical trajectory for different values of $\alpha$ and $\beta.$} 
\end{center}
\end{figure}
\end{example}

We would like to note here that the functions $U_1(x,y),$ $S_1(x,y),$ $U_2(x,y)$ and $S_2(x,y),$ obtained in this example are not the same as those given in \cite{kul}. However, they are some certain constant multiples of each other, and hence, the manifolds provided here and given in \cite{kul} are the same. So, we recover the results given in \cite{kul} by taking $\beta=0.$

\begin{example}
As another example let us keep $\alpha$ the same as in the previous example but change $\beta.$ 

For $\alpha=0.2,$ $\beta=0.5,$ we have
\begin{align*}
U_3(x,y)&= -0.4152273992 x + 0.8491364395 - 0.2923863004 y \\
& \qquad  +0.2419777563(-0.4152273992 x - 0.3508635604 + 0.7076136995 y)^2\\
& \qquad  -0.0974600586(-0.4152273992 x - 0.3508635604 + 0.7076136995 y)^3,\\
S_3(x,y)&=-0.4152273992 x - 0.3508635604 + 0.7076136995 y \\
& \qquad  +0.1961061968(-0.4152273992 x + 0.8491364395 - 0.2923863004 y)^2\\
& \qquad  +0.09806508071(-0.4152273992 x + 0.8491364395 - 0.2923863004 y)^3.
\end{align*}
and 
\begin{align*}
U_4(x,y)&= -0.3492151478 x + 1.214293633 -0.3253924261 y \\ 
& \qquad +0.3059452562(-0.3492151478 x - 0.5857063670 + 0.6746075740y)^2 \\
& \qquad -0.1066716833(-0.3492151478 x - 0.5857063670 + 0.6746075740y)^3, \\
S_4(x,y)&= -0.3492151478 x - 0.5857063670 + 0.6746075740 y \\
& \qquad + 0.1446549340(-0.3492151478 x + 1.214293633 - 0.3253924261y)^2 \\
& \qquad + 0.0525187072(-0.3492151478 x + 1.214293633 - 0.3253924261y)^3.
\end{align*}

Figure \ref{fig2} shows the graphs of the functions $U_3(x,y)=0,$ $S_3(x,y)=0,$ $U_4(x,y)=0$ and $S_4(x,y)=0$ together with a typical trajectory. As it can be seen, the trajectory follows the unstable manifold in both cases.

\begin{figure}[!ht]
\begin{center}
\subfigure[Graphs of $U_3(x,y)=0$ (blue) and $S_3(x,y)=0$ (red) for $\alpha=0.2, \beta=0.5$]{
\resizebox*{60mm}{!}{\includegraphics{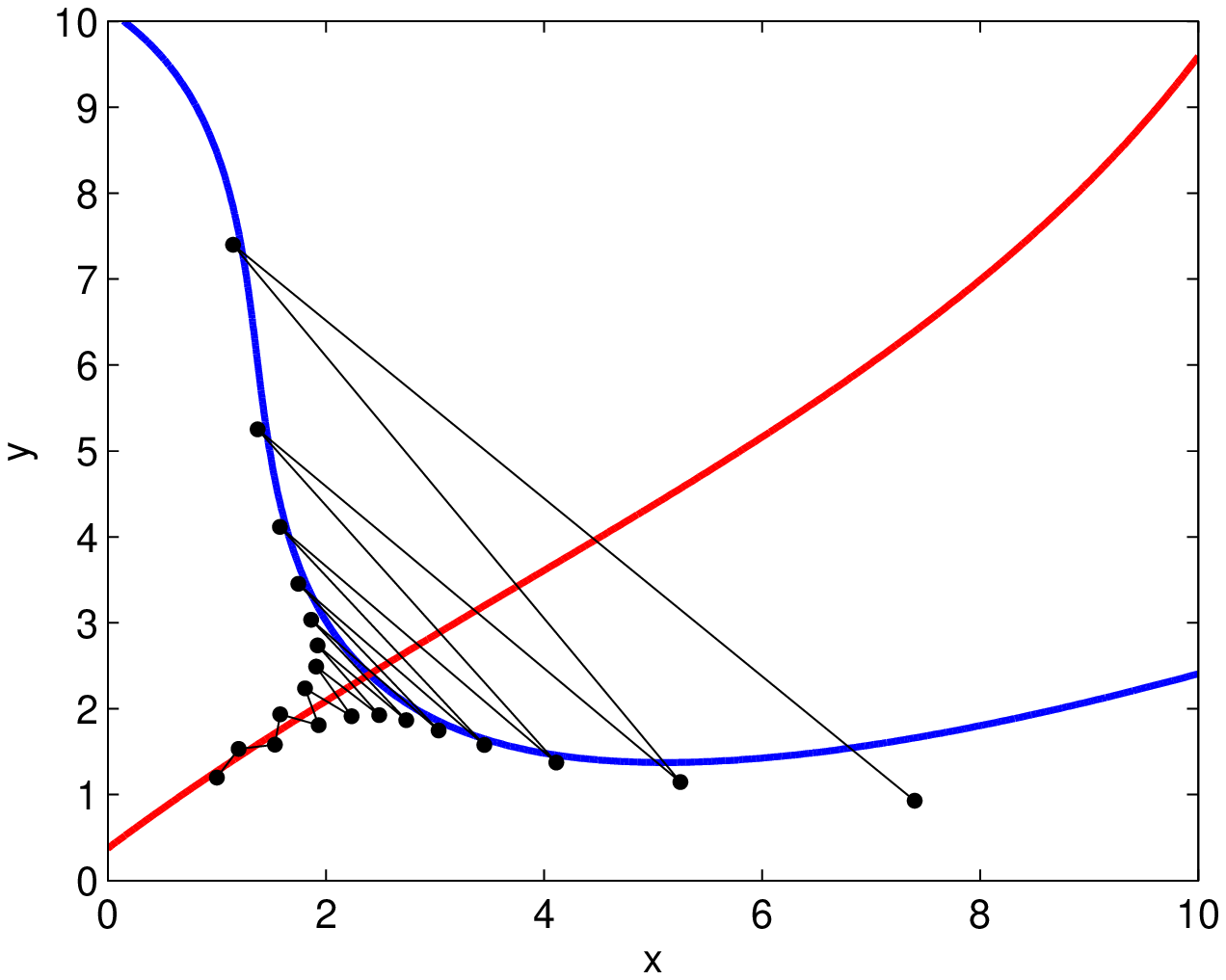}}}\hspace{8mm}%
\subfigure[Graphs of $U_4(x,y)=0$ (blue) and $S_4(x,y)=0$ (red) for $\alpha=0.8, \beta=0.5$]{
\resizebox*{60mm}{!}{\includegraphics{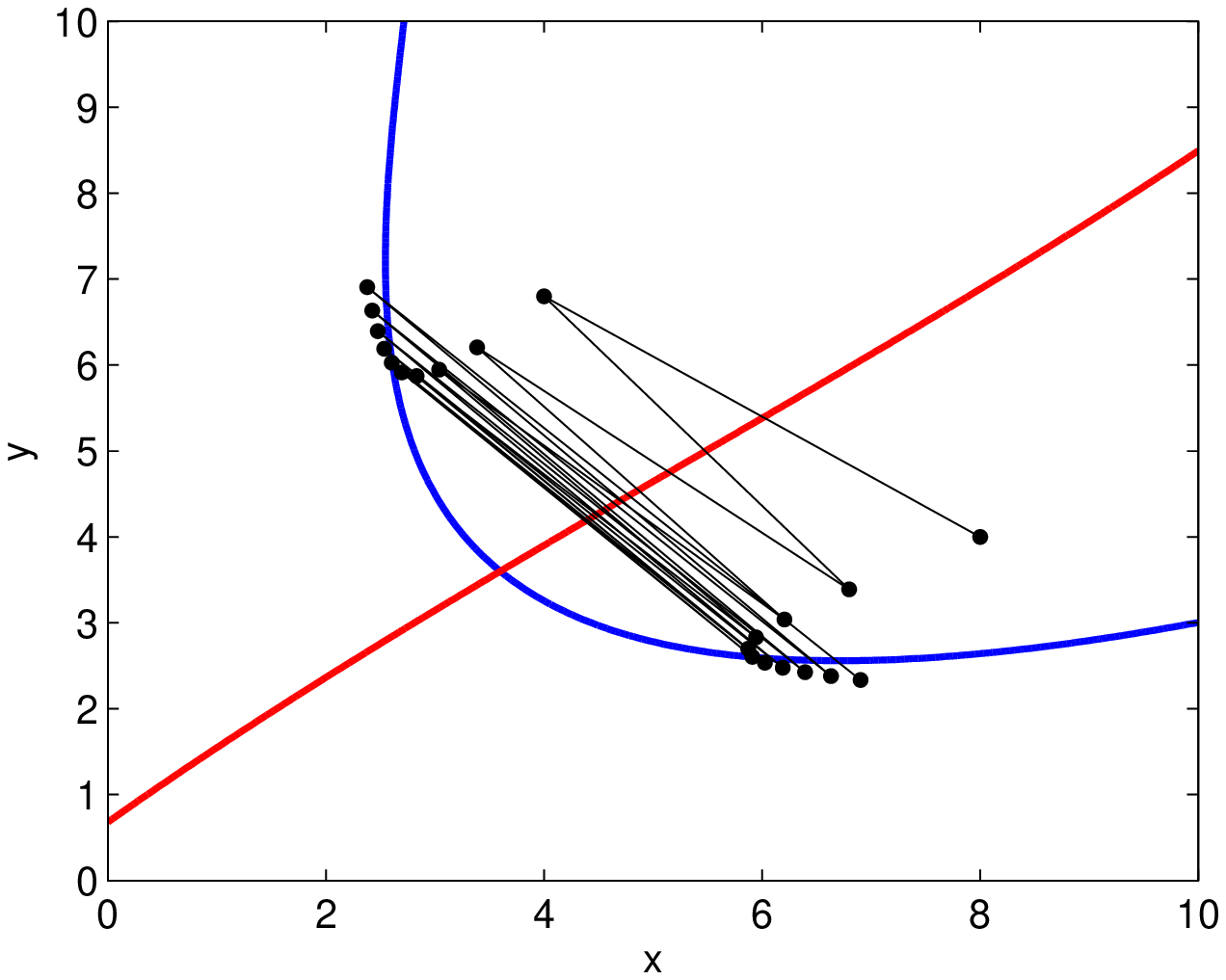}}}
\caption{\label{fig2} Graphs of stable and unstable manifolds together with a typical trajectory for different values of $\alpha$ and $\beta.$} 
\end{center}
\end{figure}
\end{example}

\begin{example} In this example, let us take the parameters as in \cite{kul}. Let $\alpha=1,$ $\beta=0.$ Let $\Phi=2.94.$ In this case, we have $\Psi=1.515463918$ and
\begin{align*}
C_1(x,y;\Phi)&=0.1491735785x+0.2260671754-0.4385703205y\\
&\qquad -0.08039102209(0.1491735785x-1.289396743+0.5614296795y)^2\\
&\qquad +0.01997063483(0.1491735785x-1.289396743+0.5614296795y)^3
\end{align*}
and
\begin{align*}
C_1(x,y;\Psi)&=0.5614296798x+1.650603257-0.8508264215y \\
&\qquad -0.1559585827(0.5614296798x-1.289396743+0.1491735785y)^2 \\ 
&\qquad +0.05514400545(0.5614296798x-1.289396743+0.1491735785y)^3.
\end{align*}
For $\Phi=2.3,$ we have $\Psi=1.769230769$ and
\begin{align*}
C_2(x,y;\Phi)&=0.2506265664x+0.4434162323-0.5764411027y\\
&\qquad -0.1137137228(0.2506265664x-1.325814536+0.4235588973y)^2\\
&\qquad +0.03453170706(0.2506265664x-1.325814536+0.4235588973y)^3
\end{align*}
and
\begin{align*}
C_2(x,y;\Psi)&=0.4235588973x+0.9741854634-0.7493734336y\\
&\qquad -0.1478278397(0.4235588973x-1.325814536+0.2506265664y)^2\\
&\qquad +0.0520650698(0.4235588973x-1.325814536+0.2506265664y)^3.
\end{align*}

Figure \ref{fig3} shows the graphs of the functions $C_1(x,y;\Phi)=0,$ $C_1(x,y;\Psi)=0$ and $C_2(x,y;\Phi)=0,$ $C_2(x,y;\Psi)=0,$ together with typical trajectories. As it can be seen, the trajectory follows the invariant manifold in both cases.
\begin{figure}[!ht]
\begin{center}
\subfigure[Graphs of $C_1(x,y;\Phi)=0$ (blue) and $C_1(x,y;\Psi)=0$ (red) for $\Phi=2.94.$]{
\resizebox*{60mm}{!}{\includegraphics{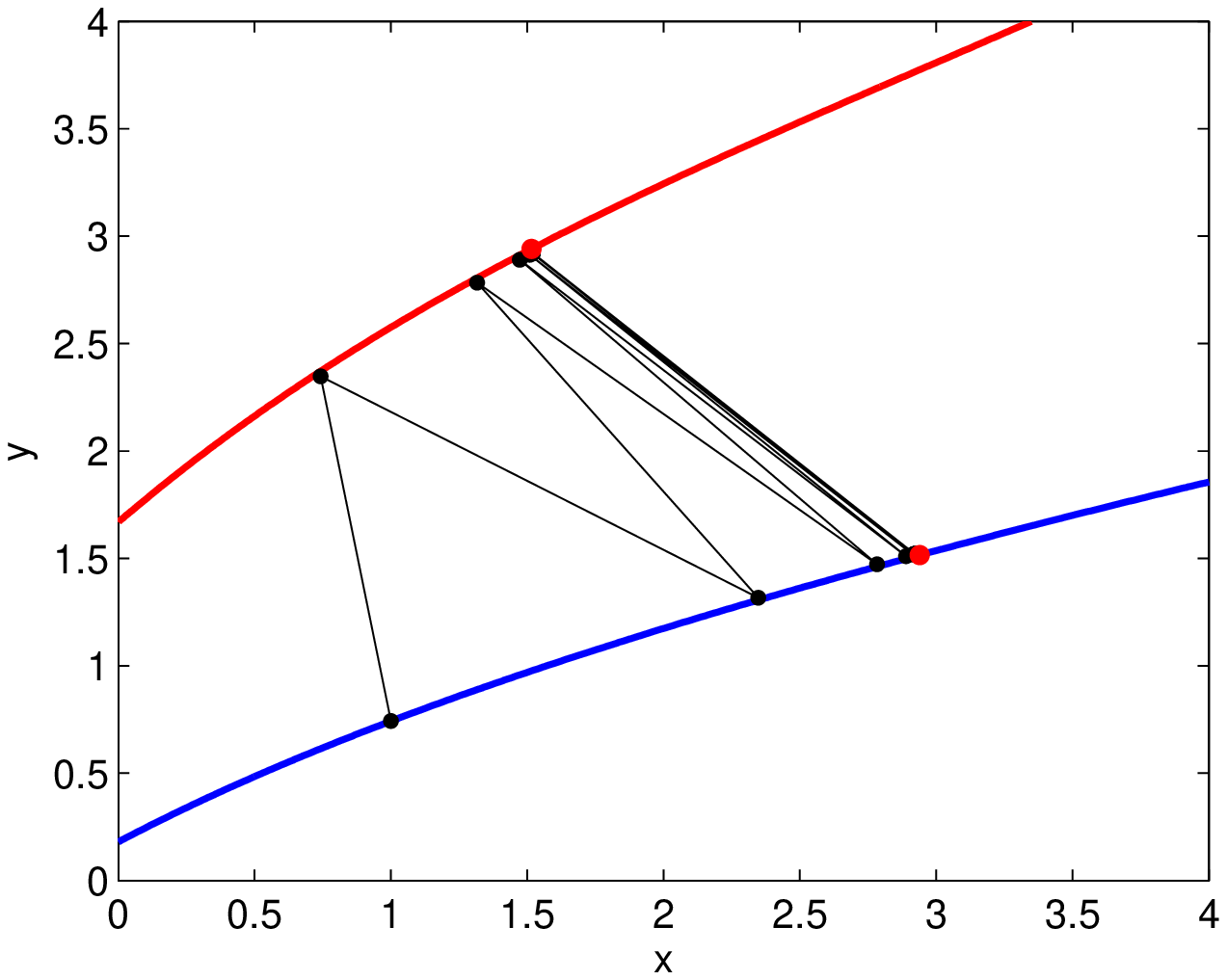}}}\hspace{8mm}%
\subfigure[Graphs of $C_2(x,y;\Phi)=0$ (blue) and $C_2(x,y;\Psi)=0$ (red) for $\Phi=2.3.$]{
\resizebox*{60mm}{!}{\includegraphics{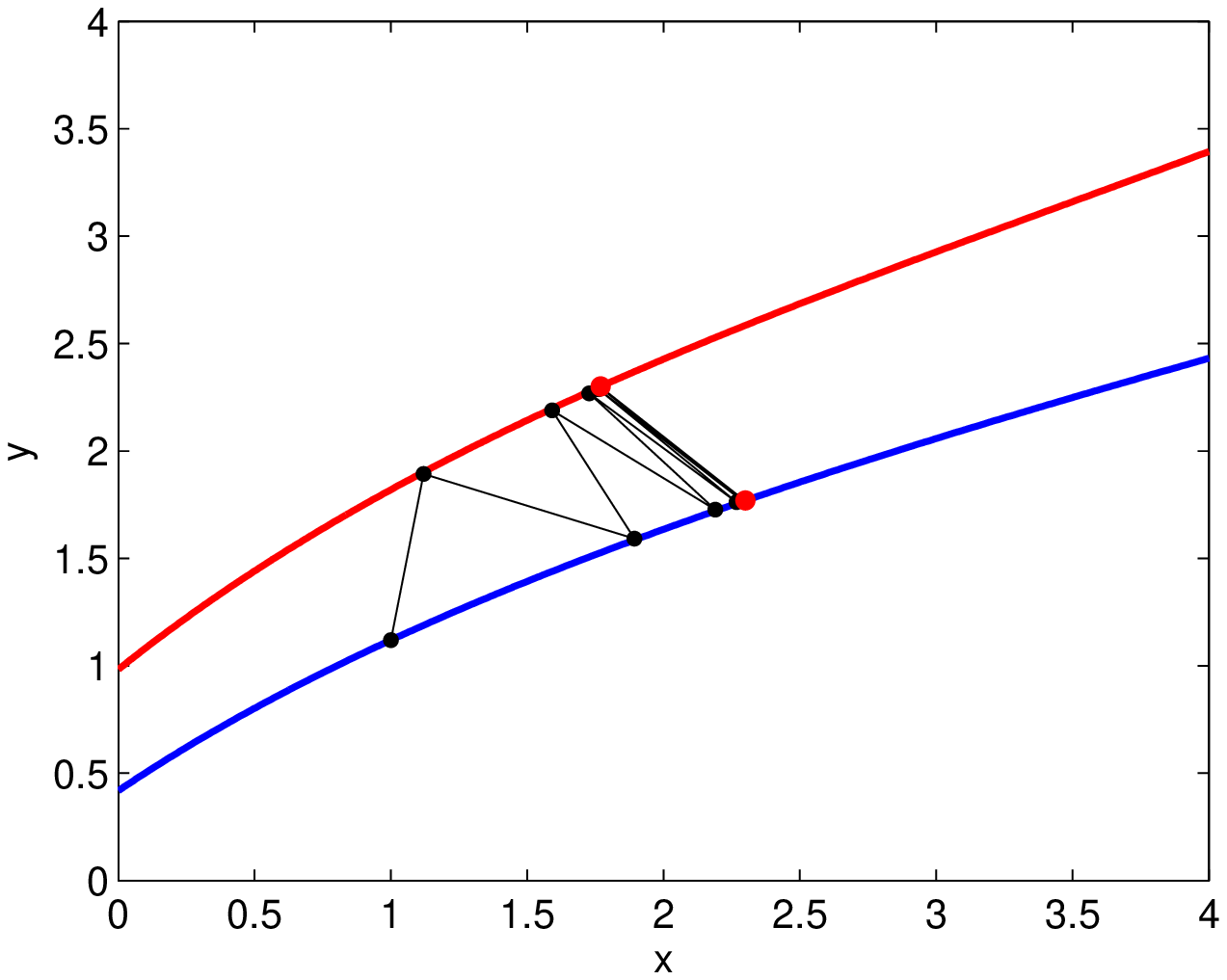}}}
\caption{\label{fig3} Graphs of invariant curves together with typical trajectories and periodic solutions for different values of $\Phi$ and $\beta=0.$} 
\end{center}
\end{figure}
\end{example}

\begin{example} As the last example, let $\alpha=1,$ $\beta=0.5.$ 
For $\Phi=2.94,$ we have $\Psi=6.255319149$ and
\begin{align*}
C_3(x,y;\Phi)&=1.071618354x+1.623998947-0.7632795057y\\
&\qquad -0.006468848599(1.071618354x-4.631320202+0.2367204943y)^2\\
&\qquad +0.001026052614(1.071618354x-4.631320202+0.2367204943y)^3
\end{align*}
and
\begin{align*}
C_3(x,y;\Psi)&=0.1416540319x+0.1686084427-0.3587413677y\\
&\qquad -0.005080796064(0.1416540319x-2.771391557+0.6412586323y)^2\\
&\qquad +0.001395071806(0.1416540319x-2.771391557+0.6412586323y)^3.
\end{align*}
For $\Phi=2.3,$ we have $\Psi=15.33333333$ and
\begin{align*}
C_4(x,y;\Phi)&=3.473613893x+6.145624586-0.9218436874y\\
&\qquad -0.001973405924(3.473613893x-9.187708748+0.07815631264y)^2\\
&\qquad +0.000140325572(3.473613893x-9.187708748+0.07815631264y)^3
\end{align*}
and
\begin{align*}
C_4(x,y;\Psi)&=0.01938877756x+0.0207414829-0.1382765531y\\
&\qquad -0.001193222187(0.01938877756x-2.279258517+0.8617234469y)^2\\
&\qquad +0.0003847285557(0.01938877756x-2.279258517+0.8617234469y)^3.
\end{align*}

Figure \ref{fig4} shows the graphs of the functions $C_3(x,y;\Phi)=0,$ $C_3(x,y;\Psi)=0$ and $C_4(x,y;\Phi)=0,$ $C_4(x,y;\Psi)=0,$ together with typical trajectories. As it can be seen, the trajectory follows the invariant manifold in both cases.
\begin{figure}[!ht]
\begin{center}
\subfigure[Graphs of $C_3(x,y;\Phi)=0$ (blue) and $C_3(x,y;\Psi)=0$ (red) for $\Phi=2.94.$]{
\resizebox*{60mm}{!}{\includegraphics{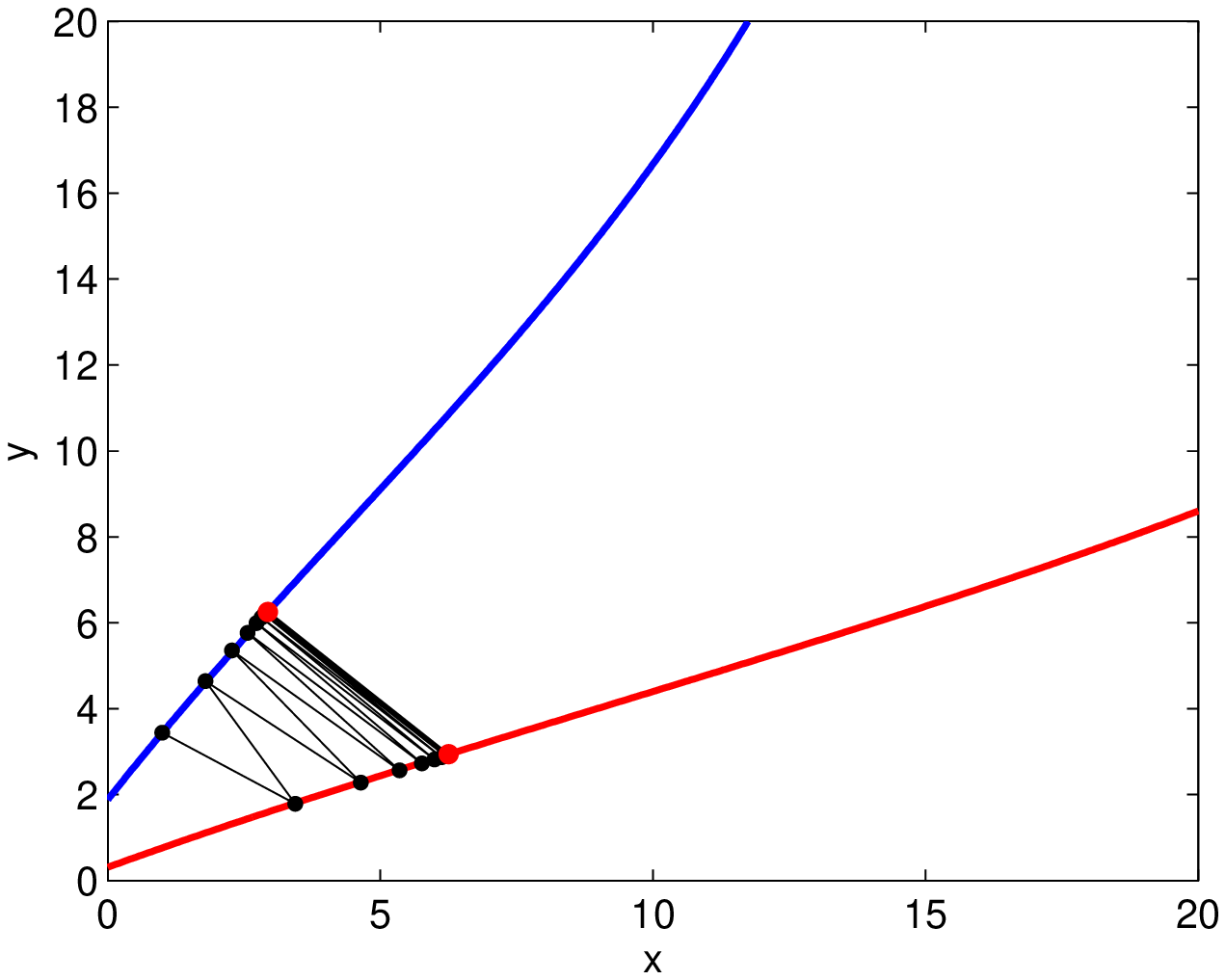}}}\hspace{8mm}%
\subfigure[Graphs of $C_4(x,y;\Phi)=0$ (blue) and $C_4(x,y;\Psi)=0$ (red) for $\Phi=2.3.$]{
\resizebox*{60mm}{!}{\includegraphics{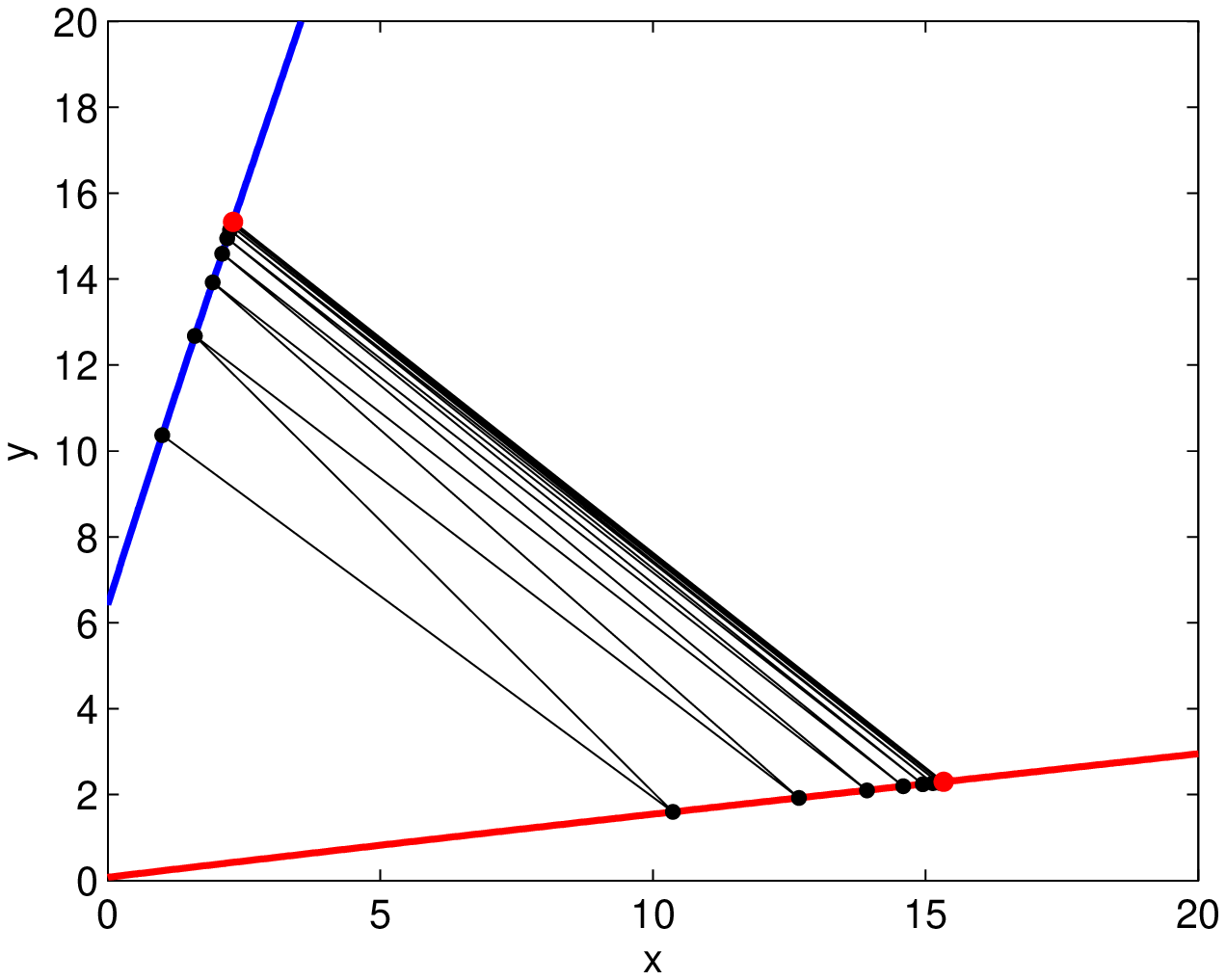}}}
\caption{\label{fig4} Graphs of invariant curves together with typical trajectories and periodic solutions for different values of $\Phi$ and $\beta=0.5.$} 
\end{center}
\end{figure}
\end{example}

%%%%%%%%%%%%%%%%%%%%%%%%%%%%%%%%%%%%%%%%%

\end{document}